\documentclass{article}

\input epsf
\usepackage[normalem]{ulem}
\usepackage{amsthm}
\usepackage{mathtools}
\usepackage{thmtools}
\declaretheoremstyle[headfont=\normalfont]{normalhead}
\usepackage{color}
\usepackage{graphicx}
\usepackage{morefloats}
\usepackage{hyperref}
\usepackage{ amssymb }
\usepackage{textcomp}
\usepackage{mathrsfs}
\usepackage{url}
\usepackage{float}
\usepackage{comment}
\allowdisplaybreaks
\usepackage{mathrsfs,bbm}

\begin{document}

\newtheorem{theorem}{Theorem}[section]
\newtheorem{corollary}[theorem]{Corollary}
\newtheorem{claim}[theorem]{Claim}
\newtheorem{conjecture}[theorem]{Conjecture}
\newtheorem{proposition}[theorem]{Proposition}
\newtheorem{lemma}[theorem]{Lemma}
\newtheorem{definition}[theorem]{Definition}
\newtheorem{question}[theorem]{Question}

\title{Linkages and removable paths avoiding vertices}

\def\qed{\hfill \rule{4pt}{7pt}}
\def\pf{\noindent {\it Proof. \/}}

\author{ \\ Xiying Du\footnote{Partially supported by a Research
    Assistantship funded by ACO and  NSF Grant DMS 1856645}, Yanjia
  Li\footnote{Supported by an ACO Research Assistantship and an ACO
    Fellowship}, Shijie Xie\footnote{Currently at Meta},  and Xingxing Yu\footnote{Partially
    supported by NSF Grant DMS 1954134}\\ \bigskip \\School of Mathematics\\Georgia Institute of Technology\\Atlanta, GA 30332
}

\maketitle

\begin{abstract}

   We say that a graph $G$ is $(2,m)$-linked if, for any distinct
  vertices $a_1,\ldots, a_m, b_1,b_2$ in $G$, there exist vertex
  disjoint connected subgraphs $A,B$ of $G$ such that $\{a_1, \ldots, a_m\}$ is
  contained in $A$ and $\{b_1,b_2\}$ is contained in $B$. A
  fundamental result in structural graph theory is the
  characterization of $(2,2)$-linked graphs, with different versions
  obtained independently by Robertson and Chakravarty, Seymour, and
  Thomassen. It appears to be very difficult to characterize
  $(2,m)$-linked graphs for $m\ge
  3$. In this paper, we provide a partial characterization of $(2,m)$-linked
  graphs by adding an average degree condition. This implies that
  $(2m+2)$-connected graphs are $(2,m)$-linked. Moreover, if $G$ is a
  $(2m+2)$-connected graph and $a_1, \ldots, a_m, b_1,b_2$ are
  distinct vertices of $G$, then there is a path $P$ in $G$ between
  $b_1$ and $b_2$ and avoiding $\{a_1, \ldots, a_m\}$ such that $G-P$
  is connected, improving a previous connectivity bound of $10m$.

\end{abstract}  

\newpage

\section{Introduction}
For any set $X$ and any integer $k\ge 0$, we let ${X\choose k}$ denote the
collection of all $k$-element subsets of $X$, and write
$\mathcal{P}(X)$ for $\bigcup_{k=0}^{|X|}{X\choose k}$, the set of all
subsets of $X$. 
For a graph $G$, we use $V(G)$ and $E(G)$ to denote its vertex and
edge set, respectively, and we use $e(G)$ and $v(G)$ to denote the number of edges and the number of vertices in $G$, respectively. For any set $A\subseteq V(G)$,  $N_G(A)$
denotes the \emph{neighborhood} of $A$ in $G$ (not including vertices of $A$), and $N_G[A]=N_G(A)\cup
A$ is the {\it closed neighborhood} of $A$ in $G$.
When there is no
confusion, we drop the subscript $G$. 
For two graphs $G$ and $H$, we use $G-H$ to denote the graph obtained
from $G$ by deleting all vertices in $V(G)\cap V(H)$ and all edges of
$G$ with at least one incident vertex in $V(G)\cap V(H)$. A path with
end vertices $u$ and $v$ is often called a $u$-$v$ {\it path}, which
is often written as a sequence of vertices (with consecutive vertices adjacent).
Let $k$ be a positive integer and $[k]:=\{1,
\ldots, k\}$. A graph $G$ is {\it $k$-linked} if for
every choice of distinct vertices $s_1, \ldots, s_k, t_1, \ldots, t_k$ of $G$, $G$
contains $k$ disjoint $s_i$-$t_i$ paths, $i\in [k]$. 

A fundamental result in graph theory is the
characterization of 2-linked graphs. 
Different versions of this result were obtained, independently, by
Robertson and Chakravarti \cite{RC79}, Seymour \cite{Se80} and
Thomassen \cite{Th80}. A linear time algorithm for finding 2-linkages
in graphs can be found in Shiloach \cite{Sh80}.  We will state Seymour's
version of the characterization, which we attempt to generalize. But
first, we need some notation. 

Let $G$ be a graph and $S\subseteq V(G)$. An {\it $S$-collection} in
$G$ is a  collection $\mathcal{X}\subseteq \mathcal{P}(V(G)\setminus
S)$ such that $N[X_1]\cap X_2=\emptyset$ for distinct
$X_1,X_2\in\mathcal{X}$.  We use $G/\mathcal{X}$ to denote the graph
obtained from $G$ by, for each $X\in\mathcal{X}$, deleting $X$ and
adding a clique on $N_G(X)$.  We say that $(G,S)$ is {\it planar} if
$G$ can be drawn in a closed disc in the plane with no edge crossings
such that the vertices in $S$ occur on the boundary of the disc. Such
a drawing is said to be a {\it disc representation} of $G$. Moreover,
if $S=\{s_1, \ldots, s_t\}$ and $s_1, \ldots, s_t$ occur on the
boundary of the disc in clockwise order then we say that $(G,s_1,
\ldots, s_t)$ is {\it planar}.

For convenience, we say that $(G,\{a_1, \ldots, a_m\}, b_1,b_2)$ is an
{\it $m$-rooted graph} or simply a {\it rooted-graph} if $G$ is a
graph and $a_1, \ldots, a_m, b_1,b_2$ are distinct vertices of $G$,
and it is said to be {\it feasible} if  $G$ contains a $b_1$-$b_2$
path $P$ such that $\{a_1,\ldots, a_m\}$ is contained in a component
of $G-P$. For convenience,  we allow $m=0$ (and hence
$\{a_1,\ldots,a_m\}=\emptyset$). 
A graph $G$ is {\it $(2,m)$-linked} if $(G,\{a_1, \ldots, a_m\},
b_1,b_2)$ is feasible for all choices of distinct vertices $a_1,
\ldots, a_m, b_1,b_2$ in $G$. Now, Seymour's characterization of
$(2,2)$-linked graphs can be restated as follows. 

\begin{theorem}[Seymour \cite{Se80}, 1980]
\label{2link} 
Let $(G, \{a_1,a_2\}, b_1,b_2)$ be a 2-rooted graph. Then either $(G, \{a_1,a_2\}, b_1,b_2)$ is feasible, or there is some $\{a_1,a_2,b_1,b_2\}$-collection $\mathcal{X}$ in $G$ such that $|N(X)|\le 3$ for all $X\in \mathcal{X}$ and $(G/\mathcal{X}, a_1,b_1,a_2,b_2)$ is planar.
\end{theorem}

In one of their papers on the graph minors project, Robertson and
Seymour \cite{RS95} generalized Theorem~\ref{2link} and 
characterized graphs $G$ with the property that  for any distinct vertices $a_1,
\ldots, a_m$, there exist $1\le i_1<i_2<i_3< i_4\le m$ such that $G$
contains disjoint
paths from $a_{i_1},a_{i_2}$ to $a_{i_3},a_{i_4}$,
respectively. It appears very difficult to find a good
characterization of feasible $m$-rooted graphs. However, if we replace
``$(G/\mathcal{X}, a_1,b_1,a_2,b_2)$ is planar'' by a bound on
$e(G/\mathcal{X})$ then we are able to obtain a more general result on
$(2,m)$-linked graphs. To state our result, we need some additional notation.

 For a rooted graph ${\cal
  G}=(G,\{a_1, \ldots, a_m\},b_1,b_2)$ and any $\{a_1, \ldots,
a_m,b_1,b_2\}$-collection ${\cal X}$ in $G$, we let ${\cal G}/{\cal X}$
denote the graph obtained from $G/{\cal X}$ by adding an edge joining
every pair of distinct vertices in $\{a_1, \ldots, a_m,b_1,b_2\}$
except the pair $\{b_1,b_2\}$.  Thus,
$$e(G/{\cal X})\le e({\cal G}/\mathcal{X})\le e(G/{\cal X})+ {m+2\choose 2}-1.$$
Our main result in this paper is the
following, which we hope would be a useful tool for potential applications.

\begin{theorem}\label{2m-link}
  Let $m\ge 0$ be an integer, and let $(G, \{a_1, \ldots  ,a_m\},b_1,b_2)$
  be an $m$-rooted graph. Then  $(G,\{a_1,
  \ldots ,a_m\},b_1,b_2)$ is feasible, or there exists some $\{a_1, \ldots ,a_m,b_1,b_2\}$-collection $\mathcal{X}$ in $G$, such that
 
  \begin{itemize}
        \item [(i)] $|N(X)|\le m+1$ for all $X\in\mathcal{X}$, and
        \item [(ii)] $e({\cal G}/\mathcal{X})\le (m+1)v(G/\mathcal{X})-m^2/2-3m/2-1$.
    \end{itemize}
    
\end{theorem}

The bound in (ii) is best possible. For integers $m\ge 0$ and $k\ge
0$, let $G_{m,k}$ be the graph on $m+2+k$ vertices $a_1,\ldots ,a_m,b_1,b_2, v_1,\ldots ,v_k$, such that $b_1v_1v_2\ldots v_kb_2$ form an induced path, each $a_i$ is adjacent to all vertices  on this path  and $G_{m,k}[\{a_1,\ldots ,a_m\}]$ consists of a $(m-1)$-clique on $\{a_1,\ldots,a_{m-1}\}$ and an isolated vertex $a_m$. 
It is easy to check that $\mathcal{G}_{m,k}=(G_{m,k},\{a_1, \ldots
,a_m\},b_1,b_2)$ is not feasible, and any $\{a_1,\ldots
,a_m,b_1,b_2\}$-collection in $G$ satisfying (i) of Theorem
\ref{2m-link} must be empty. Moreover,
\begin{align*}
    e({\cal G}_{m,k}/\emptyset)
=&e(G_{m,k})+|\{a_ia_m:i\in[m-1]\}|\\
=&(k+1)+m(k+2)+\binom{m-1}{2}+(m-1)\\ 
=&(m+1)v(G_{m,k}/\emptyset)-m^2/2-3m/2-1.
\end{align*}
So the bound in (ii) of Theorem \ref{2m-link} can be tight.

As a quick consequence of Theorem~\ref{2m-link}, we have the following

  \begin{corollary}\label{connectivity}
     Let $m\ge 0$ be an integer, and let $(G, \{a_1, \ldots
     ,a_m\},b_1,b_2)$ be an $m$-rooted graph. If $G$ is
     $(2m+2)$-connected then  $(G,\{a_1,
  \ldots ,a_m\},b_1,b_2)$ is feasible.
    \end{corollary}

Another application of Theorem~\ref{2m-link} concerns removable paths
in graphs. 
When studying 3-connected graphs and drawings of planar graphs, Tutte \cite{Tu63} proved the following result: Given a 3-connected graph $G$ and distinct vertices $a, b_1,b_2$ of $G$, there is $b_1$-$b_2$ path $P$ in $G$ such that $a\notin V(P)$ and $G-P$ is connected. A result of Jung \cite{Ju70}  implies that in a 6-connected graph $G$, for any distinct vertices $a_1, a_2,b_1,b_2$ of $G$, there is $b_1$-$b_2$ path $P$ in $G$ such that $a_1,a_2\notin V(P)$ and $G-P$ is connected. We are interested in the following more general question.

\begin{question}\label{removable}
  Let $m$ be a positive integer. Determine the least positive integer $f(m)$ with the following property:   For any $f(m)$-connected graph $G$ and  any distinct vertices $a_1, \ldots, a_m, b_1,b_2$ of $G$, there is a $b_1$-$b_2$ path $P$ in $G$ such that $\{a_1, \ldots, a_m\}\cap V(P)=\emptyset$ and $G-P$ is connected.   
\end{question}  

Tutte's result mentioned above implies $f(1)=3$, and Jung's result above implies $f(2)=6$. The third author proved $f(3)=6$ in his PhD thesis \cite{Xi19}. The exact value of $f(m)$ has not been determined for $m\ge 4$, but there are known upper bounds on $f(m)$ using results on linkages in graphs.
Let $k$ be a positive integer.  Larman and Mani \cite{LM74} and Jung \cite{Ju70} independently showed that highly connected graphs are $k$-linked. The bound on the connectivity can be significantly improved by combining results in \cite{RS95} and \cite{Ko82,Th84}, which  is further improved in \cite{BT96}. The current best bound is due to Thomas and Wollan \cite{TW05}: $10k$-connected graphs are $k$-linked. From this result, we can show that $f(m)\le 10m$, as one could use $m-1$ paths to form a connected graph containing $a_1, \ldots, a_m$ and one path to connect $b_1$ and $b_m$. 
In this paper, we improve this bound to $2m+2$. 

\begin{theorem}\label{main}
Let $m$ be any positive integer.  For any $(2m+2)$-connected graph $G$ and  any distinct vertices $a_1, \ldots, a_m, b_1,b_2$ of $G$, there is a $b_1$-$b_2$ path $P$ in $G$ such that $\{a_1, \ldots, a_m\}\cap V(P)=\emptyset$ and $G-P$ is connected. 
\end{theorem}

\medskip

We will see that to prove Theorem~\ref{main} it suffices to show that
$(G, \{a_1,\ldots,  a_m\}, b_1,b_2)$ is feasible. Our paper is
organized as follows.

In Section 2, we consider Theorem~\ref{2m-link} for
$m=0, 1$ and 2, and note that the case $m=2$ is an easy
consequence of Theorem~\ref{2link}. We also give a high level
description of the proof of Theorem~\ref{2m-link} where the concept of ``critical feasibility'' of rooted graphs arise
naturally,  and we state Theorem~\ref{criticallylinked} on critically
feasible rooted graphs.

We take care of the induction step for Theorem
\ref{2m-link} in Section 3, and the induction step of Theorem
\ref{criticallylinked} in Section 4, by assuming both Theorem \ref{2m-link} and Theorem
\ref{criticallylinked}    for smaller $m$. In Section 5, we  complete
the proofs of Theorems \ref{2m-link} and \ref{criticallylinked}. We
also prove Theorem \ref{main} in Section 5  and  offer some concluding remarks.

\section{Critical feasibility}

In this section, we give a high level discussion about
Theorem~\ref{2m-link} and its proof, where the idea of ``critical feasibility''
will arise naturally. First, we show that Theorem~\ref{2m-link} holds
for $m\le 2$. 

Let $G$ be any graph, and let $S\subseteq V(G)$.  We use $G[S]$ to
denote the subgraph of $G$ \emph{induced} by $S$ and let
$G-S=G[V(G)\setminus S]$. When $S=\{s\}$ we write $G-s$ for $G-\{s\}$. 
For any $H\subseteq G$, we write $G[H]$ for $G[V(H)]$. 

\begin{lemma}\label{base1}
    Theorem \ref{2m-link} holds for $m=0, 1, 2$. 
\end{lemma}
\pf  Let $\mathcal{G}=(G,\{a_1,\ldots,a_m\},b_1,b_2)$ be an $m$-rooted graph and suppose $\mathcal{G}$ is not feasible.

For $m=0$, $G$ contains no $b_1$-$b_2$ path. Let $D$ be the component
of $G$ containing $b_1$; so $b_2\not\in V(D)$.  Let $\mathcal{X}=\{V(D)\setminus\{b_1\}, V(G-D)\setminus \{b_2\})\}$. Then $|N(X)|\le 1$ for all $X\in\mathcal{X}$, and
$\mathcal{G}/\mathcal{X}$ has 2 vertices and no edge. Hence
$e(\mathcal{G}/\mathcal{X})=0 < 1\cdot v(G/\mathcal{X})-0^2/2-(3\cdot
0)/2-1.$

For $m=1$,  $G-a_1$ contains no $b_1$-$b_2$ path. Let $D$ be the component of $G-a_1$ containing $b_1$; then $b_2\not\in V(D)$. Let $\mathcal{X}=\{V(D)\setminus\{b_1\}, V(G-D)\setminus \{a_1,b_2\}\}$. Then $|N(X)|\le 2$ for all $X\in\mathcal{X}$, and $G/\mathcal{X}$ has 3 vertices and at most 2 edges. Hence, $e(\mathcal{G}/\mathcal{X})\le 2 <2\cdot v(G/\mathcal{X})-1^2/2-(3\cdot 1)/2-1.$

For $m=2$,  by Theorem \ref{2link}, there exists some $\{a_1,a_2,b_1,b_2\}$-collection $\mathcal{X}$ in $G$, such that $|N(X)|\le 3$ for all $X\in\mathcal{X}$, and $(G/\mathcal{X},a_1,b_1,a_2,b_2)$ is planar.
  In a disc representation of $(G/\mathcal{X},a_1,b_1,a_2,b_2)$, we can draw edges $a_1b_1,b_1a_2,a_2b_2,b_2a_1$ and $a_1a_2$ outside the disk, without introducing edge crossings.
   So  $\mathcal{G}/\mathcal{X}$ is also planar. 
   Thus, $e(\mathcal{G}/\mathcal{X})\le 3v(\mathcal{G}/\mathcal{X})-6
   =3v(G/\mathcal{X})-2^2/2-(3\cdot 2)/2-1$. \qed

   \medskip

For $m\ge 3$, we start with an $m$-rooted graph ${\cal G}:=(G, \{a_1, \ldots,
a_m\}, b_1,b_2\}$ that is not feasible. We then consider the
$(m-1)$-rooted graph ${\cal G}_{a_m}:=(G-a_m, \{a_1, \ldots,
a_{m-1}\}, b_1,b_2\}$. If ${\cal G}_{a_m}$   is not feasible then, by induction (on
$m$), there exists some $\{a_1, \ldots ,a_{m-1},b_1,b_2\}$-collection
$\mathcal{X}'$ in $G-a_m$, such that
    \begin{itemize}
        \item[(i)] $|N_{G-a_m}(X)|\le (m-1)+1$ for all $X\in\mathcal{X}'$, and
        \item[(ii)] $e(\mathcal{G}_{a_m}/\mathcal{X}')\le ((m-1)+1)v((G-a_m)/\mathcal{X}')-(m-1)^2/2-3(m-1)/2-1$.
    \end{itemize}
When viewed as a $\{a_1, \ldots, a_m,b_1,b_2\}$-collection in $G$,
one can check that ${\cal X}'$ satisfies the conclusion of
Theorem~\ref{2m-link}. Now suppose ${\cal G}_{a_m}$  is feasible and let $A, B$ be disjoint
connected subgraphs of $G-a_m$ such that $\{a_1, \ldots,
a_{m-1}\}\subseteq V(A)$ and $B$ is a $b_1$-$b_2$ path.  Hence, $\{a_1, \ldots, a_{m-1}\}$ and $a_m$ are contained in
different components of $G-B$. We choose such $B$ that the component of
$G-B$ containing $a_m$, say $D$, is maximal. 
Let $G^*=G/A$
and let $a^*$ denote
the vertex resulted from the contraction of $A$. Then $(G^*,
\{a^*,a_m\},b_1,b_2)$ is not feasible as, otherwise, we can show that
${\cal G}$ is feasible. Hence by Theorem~\ref{2link},
there is some $\{a^*,a_m,b_1,b_2\}$-collection ${\cal X}^*$ such that
$|N_{G^*}(X)|\le 3$ for all $X\in {\cal X}^*$ and $(G^*/{\cal X}^*, a^*,b_1,a_m,b_2)$ is planar.
Using this planarity, we can bound the number of edges in $G[N[D]]/\mathcal{X}_D$ for some $\mathcal{X}_D\subseteq \mathcal{X}^*$.
Let $U=N(D)\cap
V(B(b_1,b_2))$. Then for any $u\in U$, $((G-D)-u, \{a_1, \ldots,
a_{m-1}\}, b_1,b_2)$ is not feasible. This motivates the concept of
``critically feasible'' feasible graphs.

Let $m\ge 0$ and $\mathcal{G}=(G,\{a_1, \ldots ,a_m\},b_1,b_2)$.  An
ordered pair of disjoint connected subgraphs $(A,B)$ of $G$ is called
a {\it linkage pair} in $\mathcal{G}$ if $\{a_1, \ldots ,a_m\}\subseteq V(A)$ and $B$ is a $b_1$-$b_2$ path. For $m=0$, we say $(\emptyset,B)$ is a \textit{linkage pair} in $\mathcal{G}$ for any $b_1$-$b_2$ path $B\subseteq G$.
For $U\subseteq V(G)\setminus\{a_1, \ldots
,a_m,b_1,b_2\}$, $\mathcal{G}$ is said to be {\it critically feasible with respect to $U$}   if  $U\subseteq V(B)$ for every linkage pair $(A,B)$ in $\mathcal{G}$. 
Then, $\mathcal{G}$ is feasible if and only if it is critically feasible with respect to $\emptyset$.

To prove Theorem \ref{2m-link} for $m\ge 3$, we also need to prove the following statement for critically feasible rooted graphs.

\begin{theorem}\label{criticallylinked}
    Let $m\ge 0$ be an integer, let $\mathcal{G}=(G,\{a_1, \ldots ,a_{m}\},b_1,b_2)$ be an $m$-rooted graph,
    and let $U\subseteq V(G)\setminus\{a_1, \ldots ,a_{m},b_1,b_2\}$. Suppose $\mathcal{G}$
    is critically feasible with respect to $U$. 
    Then there exists some $(\{a_1, \ldots ,a_{m},b_1,b_2\}\cup U)$-collection $\mathcal{X}$ in $G$, such that
    \begin{itemize}
        \item[(i)] $|N(X)|\le m+2$ for all $X\in\mathcal{X}$, and
        \item[(ii)] $e(\mathcal{G}/\mathcal{X})\le (m+2)v(G/\mathcal{X})-m^2/2-5m/2-3-|U|$.
    \end{itemize}

  \end{theorem}

We conclude this section by proving 
Theorem~\ref{criticallylinked} for $m=0$.
 
\begin{lemma}\label{base2}
    Theorem \ref{criticallylinked} holds for all critically feasible $0$-rooted graphs.
\end{lemma}

\pf Let $\mathcal{G}=(G,\emptyset ,b_1,b_2)$ be a $0$-rooted graph
that is critically feasible with respect to a set $U\subseteq
V(G)\setminus \{b_1,b_2\}$. Let $B$ be an induced $b_1$-$b_2$ path in $G$, and
denote $U=\{u_1,\ldots,u_k\}$ such that $k=|U|$ and
$b_1,u_1,\ldots,u_k,b_2$ occur on $B$ in the order listed. Let
$u_0=b_1,u_{k+1}=b_2$. Let $$\mathcal{X}:=\{V(C): C \mbox{ is a
  component of } G-(U\cup\{b_1,b_2\})\}.$$ 

We claim that for each $X\in\mathcal{X}$, $N(X)\subseteq
\{u_i,u_{i+1}\}$ for some $i\in [k]\cup \{0\}$. 
Because otherwise there exist $u_i,u_j\in N(X)$ with $i+1<j$. 
By replacing $B[u_i,u_j]$  in $B$ with a $u_i$-$u_j$ path in $G[X\cup
\{u_i,u_j\}]$, we obtain from $B$ a $b_1$-$b_2$ path in $G-u_{i+1}$,
contradicting the assumption that $\mathcal{G}$ is critically feasible with respect to $U$.

Hence, $\mathcal{X}$ is a $(\{b_1,b_2\}\cup U)$-collection in $G$, 
$|N(X)|\le 2$ for all $X\in\mathcal{X}$, and 
$\mathcal{G}/\mathcal{X}$ is the path $b_1u_1\ldots u_kb_2$. Then \[
e(\mathcal{G}/\mathcal{X})
=k+1
=2(k+2)-3-k
=2v(G/\mathcal{X})-0^2/2-(5\cdot 0)/2-3-k,
\]
so $\mathcal{X}$ is the desired collection showing that Theorem
\ref{criticallylinked} holds for ${\cal G}$ and $U$.\qed

\medskip

We will prove Theorem \ref{2m-link} and Theorem \ref{criticallylinked} by induction on $m$. The main induction step will be divided into two lemmas given in Section 3 and Section 4.

\section{Induction step for Theorem \ref{2m-link}}

Suppose  $m\ge 3$, Theorem \ref{2m-link} holds for all $(m-1)$-rooted
graphs, and Theorem \ref{criticallylinked} holds for all critically
feasible $(m-1)$-rooted graphs.

We need to show that Theorem \ref{2m-link} holds for $m$-rooted graphs.
Let  $\mathcal{G}=(G,\{a_1, \ldots ,a_m\},b_1,b_2)$ be an $m$-rooted graph.
We apply induction on $v(G)$. Note that $v(G)\ge m+2$.
If $v(G)=m+2$, then $\mathcal{G}/\emptyset $ has $m+2$ vertices and at most $\binom{m+2}{2}=m^2/2+3m/2+1$ edges. Hence, \[
e(\mathcal{G}/\emptyset)\le m^2/2+3m/2+1=(m+1)v(G/\emptyset)-m^2/2-3m/2-1;
\]
so Theorem \ref{2m-link} holds for $\mathcal{G}$ with ${\cal X}=\emptyset$.

Now assume $v(G)\ge m+3$, and 
suppose that Theorem \ref{2m-link} holds for all $m$-rooted graphs with fewer vertices. 
Moreover, suppose, for a contradiction, that
\begin{itemize}
  \item [(1)] ${\cal G}$ is not feasible, and there does not exist any
    $\{a_1,\ldots,a_m,b_1,b_2\}$-collection $\mathcal{X}$ in $G$ satisfying (i) and (ii) of Theorem \ref{2m-link}.
 \end{itemize}

We claim that
\begin{itemize}
    \item[(2)] $|N(X)|\ge 4$ for any nonempty set $X\subseteq V(G)\setminus\{a_1,\ldots,a_m,b_1,b_2\} $.
\end{itemize}
Suppose on the contrary that there exists $X_1\subseteq V(G)\setminus\{ a_1,\ldots,a_m,b_1,b_2\} $ such that $X_1\ne\emptyset$ and $|N(X_1)|\le 3$.
We may assume $G[X_1]$ is connected; otherwise, we consider the vertex set of a component of $G[X_1]$ instead.
Let $G'=G/\{ X_1\}$,  obtained from $G$ by deleting $X_1$ and adding a
clique on $N(X_1)$, and let $\mathcal{G}'=(G', \{a_1, \ldots
,a_m\},b_1,b_2)$.

If ${\cal G}'$ is feasible then let $(A',B')$ be a linkage pair in
${\cal G}'$. Then at most one of $A',B'$  contains some edge in
$G'[N(X_1)]$, since $A'$ and $B'$ are disjoint and $|N(X_1)|\le 3$.
Then $(A',B')$, $(G[A'\cup X_1],B')$ or $(A',G[B'\cup X_1])$  is a linkage
pair in $\mathcal{G}$, showing that ${\cal G}$ is feasible,
contradicting (1).

So $\mathcal{G}'$ is not feasible. Because
$v(G')<v(G)$, we may apply the induction hypothesis and obtain some
$\{a_1,\ldots,a_m,b_1,b_2\}$-collection $\mathcal{X}'$ in $G'$,
such that $|N_{G'}(X)|\le m+1$ for all $X\in\mathcal{X}'$ and \[
e(\mathcal{G}'/\mathcal{X}')\le (m+1)v(G'/\mathcal{X}')-m^2/2-3m/2-1.
\]

If $N(X_1)\cap X=\emptyset$ for all $X\in\mathcal{X}'$, then $\mathcal{X}:=\mathcal{X}'\cup\{X_1\}$ is a $\{a_1,\ldots,a_m,b_1,b_2\}$-collection in $G$, $|N(X)|\le m+1$ for all $X\in\mathcal{X}$, and
$\mathcal{G}/\mathcal{X}=\mathcal{G}'/\mathcal{X}';$ so $
e(\mathcal{G}/\mathcal{X})\le (m+1)v(G/\mathcal{X})-m^2/2-3m/2-1.$
This contradicts (1).

Otherwise, there exist $X_2\in\mathcal{X}'$ and $v\in N(X_1)\cap
X_2$. Then, since $G'[N(X_1)]$ is a clique, $N(X_1)\setminus
\{v\}\subseteq N_{G'}[X_2]$; so $N(X_1)\cap X=\emptyset$ for $X\in
{\cal X}'\setminus \{X_2\}$ and $N(X_1\cup X_2)=N_{G'}(X_2)$. 
Let $\mathcal{X}:=(\mathcal{X}'\setminus\{X_2\})\cup\{X_1\cup X_2\}$. 
Then $\mathcal{X}$ is a $\{a_1,\ldots,a_m,b_1,b_2\}$-collection in $G$, $|N(X)|\le m+1$ for all $X\in\mathcal{X}$ and $
\mathcal{G}/\mathcal{X}=\mathcal{G}'/\mathcal{X}'$,
 which implies $
e(\mathcal{G}/\mathcal{X})\le (m+1)v(G/\mathcal{X})-m^2/2-3m/2-1
$, contradicting (1).
\qed

\medskip

We also claim that
\begin{itemize}
  \item [(3)]  ${\cal G}_{a_m}:=(G-a_m,\{a_1, \ldots ,a_{m-1}\},b_1,b_2)$ is feasible.
\end{itemize}  
For, suppose $\mathcal{G}_{a_m}$ is not feasible. By
assumption,  Theorem~\ref{2m-link} holds for ${\cal G}_{a_m}$; so
there exists some $\{a_1, \ldots ,a_{m-1},b_1,b_2\}$-collection $\mathcal{X}$ in $G-a_m$, such that 
$|N_{G-a_m}(X)|\le m$ for all $X\in\mathcal{X}$, 
and  $$e(\mathcal{G}_{a_m}/\mathcal{X})\le m\cdot v((G-a_m)/\mathcal{X})-m^2/2-m/2.$$
Now $\mathcal{X}$ is a $\{a_1, \ldots, a_m, b_1,b_2\}$-collection in $G$, 
$|N_G(X)|\le |N_{G-a_m}(X)|+1\le m+1 \mbox{ for all } X\in\mathcal{X}$,
and
\begin{align*}
    e(\mathcal{G}/\mathcal{X})&\le e(\mathcal{G}_{a_m}/\mathcal{X})+ d_{\mathcal{G}/\mathcal{X}}(a_m)\\
&\le m\cdot (v(G/\mathcal{X})-1)-m^2/2-m/2 +(v(G/\mathcal{X})-1)\\
&=(m+1)v(G/\mathcal{X})-m^2/2-3m/2-1.
\end{align*}
Therefore, $G$ and $\mathcal{X}$ satisfy (i) and (ii)  of Theorem
\ref{2m-link}, contradicting (1). \qed

\medskip

By (3),  let $(A,B)$ be a linkage pair in ${\cal G}_{a_m}$. Let $D$ denote the component of $G-B$ containing $a_m$. We choose the pair $(A,B)$ such that 
\[v(D) \mbox{ is maximal}.\]
Without loss of generality, we may assume that $B$ is induced. Denote  $G_1=G-D$ and $U=N(D)\setminus \{b_1,b_2\}=\{u_1, \ldots ,u_k\}$, such that $b_1,u_1,u_2, \ldots, u_k,b_2$ occur on $B$ in the order listed. Let $u_0=b_1$, $u_{k+1}=b_2$ and $U^+=U\cup\{b_1,b_2\}$. 
We define the graph $A_m$ by letting
$V(A_m)=V(D)\cup U^+$ and $E(A_m)=\{e\in E(G): e\cap D\ne\emptyset\}.$
Then
\begin{itemize}
  \item [(4)] $V(A_m)\cap V(G_1)= U^+$, $E(A_m)\cap E(G_1)=\emptyset$, and $E(A_m)\cup E(G_1)=E(G)$.
\end{itemize}

Let $G^*$ be the graph obtained from $G$ by contracting $A$ to a vertex $a^*$. 
If $(G^*,\{a^*,a_m\}, b_1,b_2)$ is feasible, i.e., $G^*$ contains
disjoint paths $A^*,B^*$ from $a^*,b_1$ to $a_m,b_2$, respectively,
then $(G[(A^*-a^*)\cup A], B^*)$ is a linkage pair in ${\cal G}$,
contradicting (1). So $(G^*,\{a^*,a_m\}, b_1,b_2)$ is not
feasible and, by Theorem \ref{2link}, 
there exists some $\{a^*,a_m,b_1,b_2\}$-collection $\mathcal{X}^*$ in $G^*$, such that
\begin{itemize}
\item [(5)]   $|N_{G^*}(X)|\le 3$  for all $X\in \mathcal{X}^*$, and $(G^*/\mathcal{X}^*,a_m,b_1,a^*,b_2)$ is planar.
\end{itemize}     

By choosing such $\mathcal{X}^*$ maximizing $|\mathcal{X}^*|$, we may
assume that $G[X]$ is connected for  all $X\in\mathcal{X}^*$.
By (2), $a^*\in N_{G^*}(X)$ for all $X\in \mathcal{X}^*$.
Let $\mathcal{X}_B=\{X\in\mathcal{X}^*:X\cap V(B)\ne\emptyset\}$. We claim that
\begin{itemize}
    \item[(6)]  for $X\in\mathcal{X}_B$, $N_{G^*}(X)=\{a^*,l(X),r(X)\}$ for some $l(X),r(X)\in V(B)$, $B\cap G^*[N_{G^*}[X]]=B[l(X),r(X)]$, and $X\cap V(A_m)=\emptyset$.
\end{itemize}
To prove (6), let  $b'\in X\cap V(B)$. Because  $b'\in X$ and $b_1,b_2\not\in X$, $N_{G^*}(X)$ must contain a vertex $l(X)$ in $B[b_1,b')$ and a vertex $r(X)$ in $B(b',b_2]$; so $N_{G^*}(X)=\{a^*,l(X),r(X)\}$. 
Since $a^*\not\in V(B)$ and $B$ is a path, $B\cap G^*[N_{G^*}[X]]=B[l(X),r(X)]$.

Now suppose there exists some vertex $v\in V(A_m)\cap X$. 
Note that $v\not\in\{b_1,b_2\}$.
Since $G^*[D\cup\{v\}]$ is connected and $a_m\not\in X$,
 $N_{G^*}(X)$ must contain a vertex of $D$. However, $\{a^*,l(X),r(X)\}\cap V(D)=\emptyset$,
a contradiction. 
\qed 

\medskip

By (6), $\{B(l(X),r(X)):X\in\mathcal{X}_B\}$ are pairwise disjoint subpaths of $B$ avoiding $U$.
Let $B^*$ be the path obtained from $B$ by, for each $X\in\mathcal{X}_B$, removing vertices in $B(l(X),r(X))$ and adding an edge $l(X)r(X)$.
Then
\begin{itemize}
    \item[(7)]$B^*$ is a $b_1$-$b_2$ path in $G^*/\mathcal{X}^*$, and  $b_1,u_1,\ldots,u_{k},b_2$ occur on $B^*$ in the order listed.
\end{itemize}

Now we show that
\begin{itemize}
    \item[(8)] $(A_m,b_1,u_1,u_2,\ldots,u_{k},b_2,a_m)$ is planar.
\end{itemize}
To prove (8), we first show $A_m\subseteq G^*/\mathcal{X}^*$.
Note that $A_m\subseteq G^*$ since $A_m\subseteq G$ and $ V(A_m)\cap
V(A)=\emptyset$ in $G$. Hence, if $A_m\not\subseteq
G^*/\mathcal{X}^*$ then there exist $X\in\mathcal{X}^*$ and $v\in X$
such that $v\in V(A_m)=V(D)\cup U^+$. By (6), $X\not\in\mathcal{X}_B$.
So $v\in V(D) $.
Since $X$ is connected and $a^*\in N_{G^*}(X)$, $G^*[X\cup\{a^*\}]$ contains a $a^*$-$v$ path $P$.
Recall that $a^*$ and $v\in V(D)$ are in different components of $G^*-B$,
so $P(a^*,v)$ must contain a vertex $b'\in V(B)$. This forces $X\in\mathcal{X}_B$, a contradiction.

Now we fix a disc representation of
$(G^*/\mathcal{X}^*,a_m,b_1,a^*,b_2)$ in the plane, and let $K$ be the
compact region in the plane bounded by $B^*$ together with the
boundary arc of the disc from $b_2$ to $b_1$ that contains $a_m$. Then $b_1,u_1,u_2,\ldots,u_{k},b_2,a_m$ occur on the boundary of $K$ in clockwise order.

It suffices to show  that $A_m\subseteq G^*/\mathcal{X}^*$ is drawn in the region $K$. Vertices in $U^+$ are on the boundary of $K$. For any $v\in V(A_m)\setminus U^+=V(D)$, the connected graph $(G^*/\mathcal{X}^*)[V(D)]=D$ contains a  $a_m$-$v$ path, which does not intersect $B^*$ and hence must lie inside $K$. So all vertices in $V(A_m)$ are drawn in $K$. Finally, 
 for any edge $e=uv\in E(A_m)$, by definition of $A_m$, at most one of $u,v$ is on $B^*$. Therefore $e$ is also drawn in $K$.
\qed 

\medskip

Let $A_m^+$ be the graph obtained from $A_m$ by adding edges $b_1a_m,b_2a_m$ (if they are not already present),
so that $(A_m^+,b_1,u_1,u_2,\ldots,u_{k},b_2,a_m)$ is planar. 
We can further add $2k+1$ additional  edges 
\[\{u_iu_{i+1}:0\le i\le k\}\cup
\{u_0u_i:2\le i\le k+1\}
\]
to $A_m^+$ without introducing edge crossings. 
Hence, by Euler's formula, we have
\begin{itemize}
    \item[(9)] $e(A_m^+)\le 3v(A_m)-6-2k-1.$
\end{itemize}

\medskip

Now consider $G_1$. We claim that

\begin{itemize}
    \item[(10)] ${\cal G}_1:=(G_1,\{a_1, \ldots ,a_{m-1}\},b_1,b_2)$ is critically feasible with respect to $U$.
\end{itemize}
To prove (10), note that  ${\cal G}_1$ is feasible by the choice of
$A,B$. If for any linkage pair $(A_1,B_1)$ in ${\cal G}_1$ we have
$U\subseteq V(B_1)$ then (10) holds. So assume that $(A_1,B_1)$ is a
linkage pair in ${\cal G}_1$ such that $u_i\in U\setminus V(B_1)$ for
some $i\in [k]$. If $u_i\in V(A_1)$, then $(G[A_1\cup D],
B_1)$ is a linkage pair in ${\cal G}$,   contradicting (1). So
$u_i\not\in V(A_1)$. Then the component of $G-B_1$ containing $a_m$
contains $D$ and $u_i$, contradicting our choice of $(A,B)$ that
$v(D)$ is maximal.
\qed

\medskip

By assumption, Theorem \ref{criticallylinked} holds for
$\mathcal{G}_1$ and $U$; so by (10) there is a $(\{a_1, \ldots ,a_{m-1},b_1,b_2\}\cup U)$-collection $\mathcal{X}$ in $G_1$ such that $|N_{G_1}(X)|\le m+1$ for all $X\in\mathcal{X}$, and
\begin{align*}
        e(\mathcal{G}_1/\mathcal{X})&\le (m+1)v(G_1/\mathcal{X})-m^2/2-3m/2-1-k.
    \end{align*} 
Note that $\mathcal{X}$ is a $\{a_1,\ldots,a_m,b_1,b_2\}$-collection in $G$, $|N_G(X)|= |N_{G_1}(X)|\le m+1$ for all $X\in\mathcal{X}$,
\[V(G/\mathcal{X})=V(G_1/\mathcal{X})\cup V(A_m),\]
and
\[
E(\mathcal{G}/\mathcal{X})= E(\mathcal{G}_1/\mathcal{X})\cup E(A_m^+)\cup\{a_ia_m:i\in[m-1]\}.
\]
Hence, 
$v(G/\mathcal{X})
    = v(G_1/\mathcal{X})+ v(A_m)-k-2,$
and 
\begin{align*}
    e(\mathcal{G}/\mathcal{X})
    &\le  e(\mathcal{G}_1/\mathcal{X})+e(A_m^+)+(m-1)\\
    &\le (m+1)v(G_1/\mathcal{X})-m^2/2-3m/2-1-k +(3v(A_m)-6-2k-1)+(m-1)\\
    &= (m+1)(v(G_1/\mathcal{X})+ v(A_m)-k-2 )
    +(m+1)(k+2)-(m-2)v(A_m)\\
    &\quad -m^2/2-m/2-3k-9\\
    &\le  (m+1)v(G/\mathcal{X}) 
    +(m+1)(k+2)-(m-2)(k+3)-m^2/2-m/2-3k-9\\
    &=(m+1)v(G/\mathcal{X})-m^2/2-3m/2-1.
\end{align*}
So $\mathcal{X}$ satisfies the assertions of Theorem \ref{2m-link}.

\medskip

We summarize the above discussion as a lemma.

\begin{lemma}\label{2m}
    Let $m\ge 3$ be an integer. Suppose Theorem \ref{2m-link} holds for all $(m-1)$-rooted graphs, and Theorem \ref{criticallylinked} holds for all critically feasible $(m-1)$-rooted graphs. Then Theorem \ref{2m-link} holds for all $m$-rooted graphs.
\end{lemma}

 \section{Induction step for Theorem \ref{criticallylinked}}

 Suppose $m\ge 2$ is an integer, Theorem~\ref{2m-link} holds
for $m'$-rooted graphs for any integer $m'$ with $0\le m'<m$, and
Theorem~\ref{criticallylinked} holds for critically feasible $(m-2)$-rooted graph.

Let $\mathcal{G}=(G,\{a_1, \ldots ,a_{m-1}\},b_1,b_2)$ be an
$(m-1)$-rooted graph that is critically feasible with respect to a set
$U\subseteq V(G)\setminus\{a_1, \ldots ,a_{m-1},b_1,b_2\}$.
Our aim is to find some  $(\{a_1, \ldots ,a_{m-1},b_1,b_2\}\cup U)$-collection $\mathcal{X}$ in $G$ satisfying
    \begin{itemize}
        \item[(i)] $|N(X)|\le m+1$ for all $X\in\mathcal{X}$, and
        \item[(ii)] $e(\mathcal{G}/\mathcal{X})\le (m+1)v(G/\mathcal{X})-m^2/2-3m/2-1-|U|$.
    \end{itemize}

We apply induction on $|U|$. For $|U|=0$, let $\mathcal{X}=\{X\}$
where $X=V(G)\setminus\{a_1,\ldots,a_{m-1},b_1,b_2\}$. Then $|N(X)|\le
m+1$ (as $N(X)\subseteq \{ a_1,\ldots,a_{m-1},b_1,b_2\}$), 
$v(G/\mathcal{X})=m+1$,  and $e({\cal G}/{\cal X})\le
\binom{m+1}{2}=m^2/2+m/2$. Hence,  \[
e(\mathcal{G}/\mathcal{X})\le {m+1\choose 2}=(m+1)^2-m^2/2-3m/2-1=(m+1)v(G/\mathcal{X})-m^2/2-3m/2-1-|U|;
\]
so $\mathcal{X}$ is the desired collection satisfying (i) and (ii). Thus, we may assume $|U|\ge 1$ and that

\begin{itemize}
  \item [(1)] the desired collection (i.e.,
satisfying (i) and (ii)) exists for any $(m-1)$-rooted graph
$\mathcal{G}'$ that is critically feasible with respect to any set
$U'$ with $|U'|<|U|$.
\end{itemize}

Let $(A,B)$ be a linkage pair in $\mathcal{G}$.
 Without loss of generality, we assume that $B$ is induced. 
 Denote  $U=\{u_1, \ldots ,u_k\}$ such that $b_1,u_1,u_2, \ldots ,u_k,b_2$ occur on $B$ in the order listed, and write
 $u_0=b_1$ and $u_{k+1}=b_2$. 

By definition,  the $(m-1)$-rooted graph
$\mathcal{G}_{u_1}:=(G-u_1,\{a_1, \ldots ,a_{m-1}\},b_1,b_2)$ is not
feasible. So  by assumption,  Theorem \ref{2m-link} holds for ${\cal
  G}_{u_1}$; and there exists some $\{a_1, \ldots ,a_{m-1},b_1,b_2\}$-collection $\mathcal{X}_1$ in $G-u_1$ such that
     \begin{enumerate}
         \item[(a)] $|N_{G-u_1}(X)|\le m$ for all $X\in\mathcal{X}_1$, and
         \item[(b)] $ f(\mathcal{X}_1):= e(\mathcal{G}_{u_1}/\mathcal{X}_1)-m\cdot v((G-u_1)/\mathcal{X}_1)\le -m^2/2-m/2$.
    \end{enumerate}
  We further choose $\mathcal{X}_1$ such that
  \begin{enumerate}
        \item[(c)]  subject to (a) and (b), $f(\mathcal{X}_1)$ 
          is minimized;
       \item[(d)] subject to (c), $|\bigcup_{X\in\mathcal{X}_1}X|$ is minimal; and
         \item[(e)] subject to (d), $|\mathcal{X}_1|$ is maximal.
     \end{enumerate}

     First, we observe that  
     \begin{enumerate}
       \item [(2)]   $G[X]$ is connected for each
         $X\in\mathcal{X}_1$.
       \end{enumerate}
       For, otherwise, suppose $G[X]$ is not connected for some
       $X\in {\cal X}_1$. Let ${\cal Y}$ denote the collection of the
       vertex sets of components of $G[X]$, and let ${\cal X}_1'=({\cal
         X}_1\setminus \{X\})\cup {\cal Y}$. Then $|{\cal X}_1'|>
       |{\cal X}_1|$ and $N(Y)\subseteq N(X)$ for all $Y\in {\cal
         Y}$.  Hence, it is easy to verify that
       $\mathcal{X}_1'$ satisfies (a),        $\bigcup_{X\in {\cal X}_1'}X =\bigcup_{X\in {\cal X}_1}X$, 
       $v((G-u_1)/{\cal X}_1')=v((G-u_1)/{\cal X}_1)$, and
       $E((G-u_1)/{\cal X}_1')\subseteq  E((G-u_1)/{\cal X}_1)$ (so $e({\cal
         G}_{u_1}/{\cal X}_1')\le e({\cal G}_{u_1}/{\cal X}_1)$.  Thus $f({\cal
         X}_1')\le f({\cal X}_1)$ and ${\cal X}_1'$ satisfies (b). Therefore, either $f({\cal
         X}_1')<f({\cal X}_1)$ and we have a contradiction to (c), or  $f({\cal
         X}_1')= f({\cal X}_1)$ and we derive a  contradiction to (e).\qed

       \medskip

We claim that
\begin{itemize}
    \item[(3)] for any   $X\in\mathcal{X}_1$,
      $|N_{G-u_1}(X)|=m$ and,  for any 
      $N_{G-u_1}(X)$-collection ${\cal Y} $  in $G[N_{G-u_1}[X]]$ with
      $|N_{G-u_1}(Y)|\le m$ for all $Y\in\mathcal{Y}$, 
      $(\mathcal{X}_1\setminus\{X\})\cup\mathcal{Y}$ is  an
      $\{a_1, \ldots ,a_{m-1},b_1,b_2\}$-collection in $G-u_1$
      satisfying (a) and
      $f((\mathcal{X}_1\setminus\{X\})\cup\mathcal{Y})>f(\mathcal{X})$
      when ${\cal Y}\ne \{X\}$.
\end{itemize}
To prove (3), let $X\in {\cal X}_1$. Suppose $|N_{G-u_1}(X)|=t<m$. Let $v$ be an arbitrary vertex in $X$ and let $\mathcal{X}_1'=(\mathcal{X}_1\setminus\{X\})\cup\{X\setminus\{v\}\}$. Note that $|N_{G-u_1}(X\setminus\{v\})|\le |N_{G-u_1}(X)\cup\{v\}|\le m$. Thus $\mathcal{X}_1'$ is a $\{a_1, \ldots ,a_{m-1},b_1,b_2\}$-collection in $G-u_1$ satisfying (a). Observe that \[
v((G-u_1)/\mathcal{X}_1') =v((G-u_1)/\mathcal{X}_1)+1,\quad 
e(\mathcal{G}_{u_1}/\mathcal{X}_1')\le e(\mathcal{G}_{u_1}/\mathcal{X}_1)+t.
\]
Hence,  
\begin{align*}   f(\mathcal{X}_1') 
   & = f(\mathcal{X}_1)+(e(\mathcal{G}_{u_1}/\mathcal{X}_1')-e(\mathcal{G}_{u_1}/\mathcal{X}_1))-m(v((G-u_1)/\mathcal{X}_1')-v((G-u_1)/\mathcal{X}_1))\\
   &\le  f(\mathcal{X}_1)+t-m\cdot  1\\
   & <  f(\mathcal{X}_1)\le -m^2/2-m/2.
\end{align*}
So $\mathcal{X}_1'$ satisfies (b), and $\mathcal{X}_1'$
contradicts the choice of ${\cal X}_1$ that $f(\mathcal{X}_1)$ is
minimal (see (c)). 

Now let ${\cal Y}$ be an
$N_{G-u_1}(X)$-collection in $G[N_{G-u_1}[X]]$ with  $|N_{G-u_1}(Y)|\le
m$ for all $Y\in\mathcal{Y}$. It is easy to see that
$(\mathcal{X}_1\setminus\{X\})\cup\mathcal{Y}$ is also a $\{a_1,
\ldots ,a_{m-1},b_1,b_2\}$-collection in $G-u_1$ satisfying (a).
Hence, $f((\mathcal{X}_1\setminus\{X\})\cup\mathcal{Y})\ge f(\mathcal{X}_1)$, otherwise $(\mathcal{X}_1\setminus\{X\})\cup\mathcal{Y}$ satisfies (b) as well and therefore contradicts (c). Suppose $f((\mathcal{X}_1\setminus\{X\})\cup\mathcal{Y})=f(\mathcal{X}_1)$. 
Then by (d), $|\bigcup_{X\in(\mathcal{X}_1\setminus\{X\})\cup\mathcal{Y}}X|\ge  |\bigcup_{X\in\mathcal{X}_1}X|$.
Since $\bigcup_{Y\in\mathcal{Y}}Y\subseteq X$, we must have 
$|\bigcup_{X\in(\mathcal{X}_1\setminus\{X\})\cup\mathcal{Y}}X|\le
|\bigcup_{X\in\mathcal{X}_1}X|$; so 
$\bigcup_{Y\in\mathcal{Y}}Y= X\ne\emptyset$ as members of ${\cal
  X}_1$ are pairwise disjoint.
Then by (e), $|\mathcal{Y}|\le |\{X\}|=1$, so
$\mathcal{Y}=\{X\}$. Thus, if ${\cal Y}\ne \{X\}$ then $f((\mathcal{X}_1\setminus\{X\})\cup\mathcal{Y})>f(\mathcal{X}_1)$.
\qed

\medskip

For $X\in\mathcal{X}_1$ with $N[X]\cap V(B)\ne\emptyset,$
let $l(X),r(X)$ denote the vertices in $N[X]\cap V(B)$
such that $B[b_1,l(X)]$ and $B[r(X),b_2]$ are minimal, and we say that
$X$ is {\it $U$-relevant} if both $N[X]\cap V(B)\ne \emptyset$ and $V(B(l(X),r(X)))\cap U\ne \emptyset$.
Then
\begin{itemize}
    \item[(4)] $l(X_1)=u_1$ for all $U$-relevant $X_1\in {\cal X}_1$. 
    \end{itemize}
Let $X_1\in {\cal X}_1$ be $U$-relevant and $H=G[N_{G-u_1}[X_1]]$. 
For convenience, write $l=l(X_1)$ and $r=r(X_1)$. To prove (4), suppose for
a contradiction that $l\ne u_1$. Note that $u_1\ne r$ since
$V(B(l,r))\cap U\ne\emptyset$; so $l,r\in N_{G-u_1}(X_1)\subseteq
V(H)$. Consider the $(m-2)$-rooted graph 
$\mathcal{H}=(H,N_{G-u_1}(X_1)\setminus\{l,r\},l,r)$.  

Suppose ${\cal H}$ is not feasible.  Then by assumption, Theorem
\ref{2m-link} holds for ${\cal H}$;   so  there exists an $N_{G-u_1}(X_1)$-collection $\mathcal{Y}$ in $H$ such that $|N_{G-u_1}(Y)|=|N_{H}(Y)|\le m-1$ for all $Y\in\mathcal{Y}$, and 
\[
 e(\mathcal{H}/\mathcal{Y})\le  (m-1)v(H/\mathcal{Y})-m^2/2+m/2.
\]
Let
$\mathcal{X}_1':=(\mathcal{X}_1\setminus\{X_1\})\cup\mathcal{Y}$. Then
$\mathcal{X}_1'$ is a $\{a_1, \ldots ,a_{m-1},b_1,b_2\}$-collection
in $G-u_1$ satisfying (a). Since  $|N_{H}(Y)|< m=|N_{H}(X_1)|$ for
$Y\in\mathcal{Y}$, we have $\mathcal{Y}\ne \{X_1\}$. Then by (3),
$f(\mathcal{X}_1')>f(\mathcal{X}_1).$ Since $G[X_1]$ is connected (by
(2)), $v(H/\mathcal{Y})\ge m+1$. Note that
$\mathcal{G}_{u_1}/\mathcal{X}_1'\subseteq
\mathcal{G}_{u_1}/\mathcal{X}_1\cup \mathcal{H}/\mathcal{Y}$, and
$E(\mathcal{G}_{u_1}/\mathcal{X}_1')\cap
E(\mathcal{H}/\mathcal{Y}) \supseteq \binom{N_{G-u_1}(X_1)}{2}\setminus\{lr\}$.
Hence,
\[
e(\mathcal{G}_{u_1}/\mathcal{X}_1')\le e(\mathcal{G}_{u_1}/\mathcal{X}_1)
+e(\mathcal{H}/\mathcal{Y})- \left(\binom{m}{2}-1\right)
=e(\mathcal{G}_{u_1}/\mathcal{X}_1)
+e(\mathcal{H}/\mathcal{Y})-m^2/2+m/2+1.
\]
Thus,
\begin{align*}
    f(\mathcal{X}_1') - f(\mathcal{X}_1)
    &=(e(\mathcal{G}_{u_1}/\mathcal{X}_1')-e(\mathcal{G}_{u_1}/\mathcal{X}_1))-m(v(\mathcal{G}_{u_1}/\mathcal{X}_1')-v(\mathcal{G}_{u_1}/\mathcal{X}_1))\\
    &\le e({\cal H}/{\cal Y})-m^2/2+m/2+1 -m(v(H/\mathcal{Y})-m )\\
    &\le (m-1)v(H/\mathcal{Y})-m^2/2+m/2-m^2/2+m/2+1 -m(v(H/\mathcal{Y})-m )\\
    &= -v(H/\mathcal{Y})+m+1\\
    &\le 0,
\end{align*}
contradicting (3).

Therefore, $\mathcal{H}$ is feasible. Let $(A_\mathcal{H},B_\mathcal{H})$ be a linkage pair in $\mathcal{H}$. Then
\begin{equation*}
    A':=\begin{cases}
A&\textrm{ if $V(A)\cap N(X_1)=\emptyset$}\\
G[(V(A-H)) \cup V(A_\mathcal{H})]& \textrm{ if $V(A)\cap N(X_1)\ne\emptyset$ }
\end{cases},
\end{equation*}
and
\begin{equation*}
  B':=B[b_1,l]\cup B_\mathcal{H}\cup B[r,b_2]  
\end{equation*}
form a linkage pair in $\mathcal{G}$. Since $\mathcal{G}$ is
critically feasible with respect to $U$, we have $U\subseteq
V(B')=V(B[b_1,l])\cup V(B_\mathcal{H})\cup V(B[r,b_2])$. So
$V(B(l,r))\cap U\subseteq V(B_\mathcal{H})$. Therefore, $\mathcal{H}$
is critically feasible with respect to $U':=V(B(l,r))\cap U$.
By induction hypothesis,  Theorem \ref{criticallylinked} holds for
$\mathcal{H}$ and $U'$; hence there exists some $(N_{G-u_1}(X_1)\cup U')$-collection $\mathcal{Y}'$ in $H$ such that $|N_{G-u_1}(Y)|=|N_{H}(Y)|\le m$ for all $Y\in\mathcal{Y}'$, and
 $$e(\mathcal{H}/\mathcal{Y}')\le m\cdot
 v(H/\mathcal{Y}')-m^2/2-m/2-|U'|.$$
 
Let $\mathcal{X}_1'':=(\mathcal{X}_1\setminus\{X_1\})\cup\mathcal{Y}'$. Then $\mathcal{X}_1''$ is also a $\{a_1, \ldots ,a_{m-1},b_1,b_2\}$-collection in $G-u_1$ satisfying (a). 
Note that $\mathcal{Y}'\ne\{X_1\}$, because  $U'\subseteq X_1$ and
each member of $\mathcal{Y}'$ is disjoint from $U'$. So by (3), $f(\mathcal{X}_1'')>f(\mathcal{X}_1).$
Further note that  $\mathcal{G}_{u_1}/\mathcal{X}_1''\subseteq
\mathcal{G}_{u_1}/\mathcal{X}_1\cup \mathcal{H}/\mathcal{Y}'$ and $E(\mathcal{G}_{u_1}/\mathcal{X}_1'')\cap
E(\mathcal{H}/\mathcal{Y}')\supseteq
\binom{N_{G-u_1}(X_1)}{2}\setminus\{lr\}$. Hence, we have
\[
e(\mathcal{G}_{u_1}/\mathcal{X}_1'')\le e(\mathcal{G}_{u_1}/\mathcal{X}_1)
+e(\mathcal{H}/\mathcal{Y}')- \left(\binom{m}{2}-1\right)
=e(\mathcal{G}_{u_1}/\mathcal{X}_1)
+e(\mathcal{H}/\mathcal{Y}')-m^2/2+m/2+1.
\]
Therefore, 
\begin{align*}
    f(\mathcal{X}_1'') - f(\mathcal{X}_1)
    &=(e(\mathcal{G}_{u_1}/\mathcal{X}_1'')-e(\mathcal{G}_{u_1}/\mathcal{X}_1))-m(v(\mathcal{G}_{u_1}/\mathcal{X}_1'')-v(\mathcal{G}_{u_1}/\mathcal{X}_1))\\
    &\le e(\mathcal{H}/\mathcal{Y}')-m^2/2+m/2+1 -m(v(H/\mathcal{Y}')-m )\\
    &\le m\cdot v(H/\mathcal{Y}')-m^2/2-m/2-|U'|-m^2/2+m/2+1 -m(v(H/\mathcal{Y}')-m)\\
    &= -|U'|+1\\
    &\le 0,
\end{align*}
a contradiction.
\qed 

\medskip

We further claim that
\begin{itemize}
    \item [(5)] for $U$-relevant $X_1\in {\cal X}_1$,
      $V(B(u_1,r(X_1)))\cap U\subseteq X_1$, and
      $\Tilde{\mathcal{H}}:=(G[N_{G}[X_1]], N_{G}(X_1)\setminus\{u_1,r\},u_1,r)$ is critically feasible with respect to $V(B(u_1,r(X_1)))\cap U$.
\end{itemize}
To prove (5), we first note that (4) implies $N(u_1)\cap
X_1\ne\emptyset$; so there exists $l'\in N(u_1)\cap X_1$. For convenience, let $r=r(X_1)$,
$H=G[N_{G-u_1}[X_1]]$ and $\tilde{H}=G[N_G[X_1]]$.
Note $H=\tilde{H}-u_1$. 
Consider the $(m-1)$-rooted graph
$\mathcal{H}':=(H,N_{G-u_1}(X)\setminus\{r\},l',r)$. 

Suppose $\mathcal{H}'$ is not feasible. Then, since Theorem
\ref{2m-link} holds for ${\cal H}'$ (by assumption),  there exists some
$(N_{G-u_1}(X)\cup\{l'\})$-collection $\mathcal{Y}$ in $H$ such that
$|N_{H}(Y)|\le m$ for all $Y\in\mathcal{Y}$ and
\begin{align*}
    e(\mathcal{H}'/\mathcal{Y})&\le m\cdot v(H/\mathcal{Y})-m^2/2-m/2.
\end{align*}
Let $\mathcal{X}_1':=(\mathcal{X}_1\setminus\{X_1\})\cup\mathcal{Y}$. Then $\mathcal{X}_1'$ is a $\{a_1, \ldots ,a_{m-1},b_1,b_2\}$-collection in $G-u_1$ satisfying (a). 
Because  $l'\in X_1$ but $l'\not\in Y$ for all $Y\in\mathcal{Y}$, we
have $\mathcal{Y}\ne \{X_1\}$. Hence,  by (3),
$f(\mathcal{X}_1')>f(\mathcal{X}_1)$. On the other hand, since
$\mathcal{G}_{u_1}/\mathcal{X}_1'\subseteq
\mathcal{G}_{u_1}/\mathcal{X}_1\cup \mathcal{H}'/\mathcal{Y}$ and
$E(\mathcal{G}_{u_1}/\mathcal{X}_1')\cap
E(\mathcal{H}'/\mathcal{Y})=\binom{N_{G-u_1}(X)}{2}$, we have
\[
e(\mathcal{G}_{u_1}/\mathcal{X}_1')\le e(\mathcal{G}_{u_1}/\mathcal{X}_1)
+e(\mathcal{H}'/\mathcal{Y})- \binom{m}{2}
=e(\mathcal{G}_{u_1}/\mathcal{X}_1)
+e(\mathcal{H}'/\mathcal{Y})-m^2/2+m/2.
\]
Therefore, 
\begin{align*}
    f(\mathcal{X}_1') - f(\mathcal{X}_1)
    &=(e(\mathcal{G}_{u_1}/\mathcal{X}_1')-e(\mathcal{G}_{u_1}/\mathcal{X}_1))-m(v(\mathcal{G}_{u_1}/\mathcal{X}_1')-v(\mathcal{G}_{u_1}/\mathcal{X}_1))\\
    &\le e(\mathcal{H}'/\mathcal{Y})-m^2/2+m/2 -m(v(H/\mathcal{Y})-m )\\
    &\le m\cdot v(H/\mathcal{Y})-m^2/2-m/2 -m^2/2+m/2 -m(v(H/\mathcal{Y})-m )\\
    &= 0,
\end{align*}
a contradiction.

So $\mathcal{H}'$ is feasible. Let
$(A_{\mathcal{H}'},B_{\mathcal{H}'})$ be a linkage pair in
$\mathcal{H}'$. Then $(A_{\mathcal{H}'},u_1l'\cup B_{\mathcal{H}'})$
is a linkage pair in $\tilde{\mathcal{H}}$;  so $\Tilde{\mathcal{H}}$ is feasible.
Moreover, for any linkage pair $(A_{\Tilde{\mathcal{H}}},B_{\Tilde{\mathcal{H}}})$ in $\Tilde{\mathcal{H}}$, 
\begin{equation*}
    A':=\begin{cases}
A&\textrm{ if $V(A)\cap N(X_1)=\emptyset$}\\
G[(V(A-\tilde{H})) \cup V(A_{\Tilde{\mathcal{H}}})]& \textrm{ if $V(A)\cap N(X_1)\ne\emptyset$ }
\end{cases},
\end{equation*}
and
\begin{equation*}
  B':=B[b_1,u_1]\cup B_{\Tilde{\mathcal{H}}}\cup B[r,b_2]  
\end{equation*}
form a linkage pair in $\mathcal{G}$. Since $\mathcal{G}$ is
critically feasible with respect to $U$, we have $U\subseteq
V(B')$. In particular, $V(B(u_1,r))\cap U\subseteq V(B_{\Tilde{\mathcal{H}}})$. Therefore, $V(B(u_1,r))\cap U \subseteq X_1$ and $\Tilde{\mathcal{H}}$ is critically feasible with respect to $V(B(u_1,r))\cap U$.
\qed 

\medskip

In the rest of the proof, let's assume Theorem \ref{criticallylinked}
does not hold for $\mathcal{G}$ and $U$. Recall the labeling of
vertices in $U$ as $u_1, \ldots, u_k$ from $b_1$ to $b_2$ in the order
listed.  We claim that
\begin{itemize}
    \item[(6)]  there exists  some $(\{a_1, \ldots ,a_{m-1},b_1,b_2\}\cup U)$-collection $\mathcal{X}_2$ in $G$ such that
      $|N(X)|\le m+1$ for all $X\in\mathcal{X}_2$,
      $e(\mathcal{G}/\mathcal{X}_2)= (m+1)v(G/\mathcal{X}_2)-m^2/2-3m/2-|U|$, and $\{a_1,\ldots,a_{m-1}\}\subseteq N_{\mathcal{G}/\mathcal{X}_2}(u_1)$.
    \end{itemize}
  We will derive ${\cal X}_2$ from ${\cal X}_1$.   First, suppose that
  no member of $\mathcal{X}_1$ is $U$-relevant.
  Then $\mathcal{X}_1$ is a $(\{a_1, \ldots ,a_{m-1},b_1,b_2\}\cup
  U)$-collection in $G$. Moreover,  $u_1$ is not adjacent to $u_i$ for any
  $i\in \{3, \ldots, k+1\}$ in $\mathcal{G}/\mathcal{X}_1$, because
  otherwise there exists a set $X$ such that $u_1,u_i\in N[X]$ and
  hence $u_2\in V(B(l(X),r(X)))\cap U$, a contradiction. Therefore, 
\begin{align*}
    e(\mathcal{G}/\mathcal{X}_1)
    &= e(\mathcal{G}_{u_1}/\mathcal{X}_1) +d_{\mathcal{G}/\mathcal{X}_1}(u_1)\\
    &\le
      mv((G-u_1)/\mathcal{X}_1)-m^2/2-m/2+(v(G/\mathcal{X}_1)-|\{u_i:
      i \in
      [k+1]\setminus \{2\}\}|\\
    &=  (m+1)v(G/\mathcal{X}_1)-m^2/2-3m/2-k.
\end{align*}
If
$e(\mathcal{G}/\mathcal{X}_1)<(m+1)v(G/\mathcal{X}_1)-m^2/2-3m/2-k$,
then $\mathcal{X}_1$ shows that Theorem
\ref{criticallylinked} holds for ${\cal G}$ and $U$, a contradiction. Therefore, the inequality
above holds with equality, which implies \[
N_{\mathcal{G}/\mathcal{X}_1}(u_1)=V(G/\mathcal{X}_1)\setminus\{u_i:
i\in [k+1]\setminus \{2\}\}.
\]
In particular, $\{a_1,\ldots,a_{m-1}\}\subseteq N_{\mathcal{G}/\mathcal{X}_1}(u_1)$. So $\mathcal{X}_2:=\mathcal{X}_1$ satisfies (6).

Thus, we may assume that there exists $X_1\in\mathcal{X}_1$ that is $U$-relevant.
Note that (4) and (5) imply that $u_1=l(X_1)$ and $u_2\in
X_1$. Since the members of $\mathcal{X}_1$ are pairwise disjoint,
$X_1$ is the unique member of $\mathcal{X}_1$ that is $U$-relevant.
As before, we write $r=r(X_1)$, $\tilde{H}:=G[N_{G}[X_1]]$, and $\Tilde{\mathcal{H}}=(\Tilde{H},N_{G}(X_1)\setminus\{u_1,r\},u_1,r)$. 
Let $U_1=U\cap V(G/\mathcal{X}_1)$ and $U_2=U\cap X_1$; then
$U=U_1\sqcup U_2$.

We claim that in the graph $\mathcal{G}/\mathcal{X}_1$, $u_1$ is not adjacent to any vertex in $(U_1\setminus\{r\})\cup\{u_{k+1}\}$. 
Otherwise, let  $u_i\in ((U_1\setminus\{r\})\cup\{u_{k+1}\})\cap N_{\mathcal{G}/\mathcal{X}_1}(u_1)$. 
Then $i\ge 3$ as $u_2\in X_1$. Since $B$ is
induced, for $u_1u_i$ to be an edge in $\mathcal{G}/\mathcal{X}_1$, $
\mathcal{X}_1$ must contain a set $X$ such that $u_1,u_i\in N(X)$.
Then $B(l(X),r(X))\cap U\supseteq B(u_1,u_i)\cap U\ne \emptyset$;  so $X=X_1$. 
Now,  since $u_i\not\in\{u_1,r\}$, it follows from (5) that $u_i\in V(B(u_1,r))\cap U\subseteq X_1$; so $u_i\not\in V(G/\mathcal{X}_1)$, a contradiction.

Hence, we have
\begin{align*}
    e(\mathcal{G}/\mathcal{X}_1)
    &= e(\mathcal{G}_{u_1}/\mathcal{X}_1) +d_{\mathcal{G}/\mathcal{X}_1}(u_1)\\
    &\le  m\cdot v((G-u_1)/\mathcal{X}_1)-m^2/2-m/2+(v(G/\mathcal{X}_1)-|(U_1\setminus\{r\})\cup\{u_{k+1}\}|)\\
    &\le    (m+1)v(G/\mathcal{X}_1)-m^2/2-3m/2-|U_1|.
\end{align*}

By assumption, Theorem \ref{criticallylinked} holds for  $\tilde{\mathcal{H}}$
and $U_2$. So by (5), there exists some $(N(X_1)\cup U_2)$-collection  $\mathcal{Y}$ in $\tilde{H}$ such that $|N(Y)|=|N_{\tilde{H}}(Y)|\le m+1$ for all $Y\in\mathcal{Y}$, and
\[e(\mathcal{\tilde{H}}/\mathcal{Y})\le (m+1)v(\tilde{H}/\mathcal{Y})-m^2/2-3m/2-1-|U_2|.
\]
Let $\mathcal{X}_2=(\mathcal{X}_1\setminus\{X_1\})\cup\mathcal{Y}$;
then $\mathcal{X}_2$ is a $(\{a_1, \ldots ,a_{m-1},b_1,b_2\}\cup U)$-collection in $G$, 
and $|N(X)|\le m+1$  for all $X\in\mathcal{X}_2$.
Note that $\mathcal{G}/\mathcal{X}_2\subseteq
\mathcal{G}/\mathcal{X}_1\cup \mathcal{\tilde{H}}/\mathcal{Y}$ and 
$\binom{N(X_1)}{2}\setminus\{u_1r\}\subseteq
E(\mathcal{G}/\mathcal{X}_1)\cap (\mathcal{\Tilde{H}}/\mathcal{Y})$.  
Further note that, 
for any $X\in\mathcal{X}_1\setminus\{X_1\}$, $N(X)$ cannot contain
both $u_1$ and $r$ (otherwise $X$ is $U$-relevant); so 
$u_1r\in E(\mathcal{G}/\mathcal{X}_2)$ if and only if $u_1r\in
E(\mathcal{\tilde{H}}/\mathcal{Y})$. Hence
\begin{align*}
& \quad e(\mathcal{G}/\mathcal{X}_2)\\
   & \le
     e(\mathcal{G}/\mathcal{X}_1)+e(\mathcal{\tilde{H}}/\mathcal{Y})-{m+1\choose
  2}\\
   &\le(m+1)v(G/\mathcal{X}_1)-m^2/2-3m/2-|U_1| 
   +(m+1)v(\tilde{H}/\mathcal{Y})-m^2/2-3m/2-1-|U_2| -{m+1\choose 2}\\
   &= (m+1)(v(G/\mathcal{X}_1)+v(\tilde{H}/\mathcal{Y})-m-1)-m^2/2-3m/2-|U_1|-|U_2|\\
&= (m+1)v(G/\mathcal{X}_2)-m^2/2-3m/2-k.
\end{align*}
If $e(\mathcal{G}/\mathcal{X}_2)<
(m+1)v(G/\mathcal{X}_2)-m^2/2-3m/2-k$, then $\mathcal{X}_2$ shows that
Theorem \ref{criticallylinked} holds for ${\cal G}$ and $U$, a contradiction.

So $e(\mathcal{G}/\mathcal{X}_2)=
(m+1)v(G/\mathcal{X}_2)-m^2/2-3m/2-k$ and, hence, all inequalities
above are in fact equalities. This means that all edges in
$\binom{N(X_1)}{2}\setminus\{u_1r\}$ belong to
$\mathcal{G}/\mathcal{X}_2$, and therefore,
$E(\mathcal{G}/\mathcal{X}_1)\setminus \{u_1r\}\subseteq
E(\mathcal{G}/\mathcal{X}_2)$. 
Since $d_{\mathcal{G}/\mathcal{X}_1}(u_1)
=v(G/\mathcal{X}_1)-|(U_1\setminus\{r\})\cup\{u_{k+1}\}|$,
 $\{a_1,\ldots,a_{m-1}\}\subseteq N_{\mathcal{G}/\mathcal{X}_1}(u_1)$.
 Then $\{a_1,\ldots,a_{m-1}\}\subseteq
 N_{\mathcal{G}/\mathcal{X}_2}(u_1)$. 
So $\mathcal{X}_2$ gives the desired collection in (6).
\qed 

\medskip
We choose the $(\{a_1, \ldots ,a_{m-1},b_1,b_2\}\cup U)$-collection
$\mathcal{X}_2$ in $G$  in  (6) such that 
\begin{itemize}
    \item [(a')] $|\bigcup_{X\in\mathcal{X}_2}X|$ is minimal; and
    \item [(b')] subject to (a'), $|\mathcal{X}_2|$ is maximal.
\end{itemize}
Then \begin{itemize}
    \item[(7)] $G[X]$ is connected for each $X\in\mathcal{X}_2$.
\end{itemize}
To see (7), suppose for a contradiction that there exists
$X_2\in\mathcal{X}_2$ such that $G[X_2]$ is not connected. Let $C_1$ be
the vertex set of some component of $G[X_2]$ and let $C_2=X_2\setminus
C_1$. Then
$\mathcal{X}_2':=(\mathcal{X}_2\setminus\{X_2\})\cup\{C_1,C_2\}$ is a
$(\{a_1, \ldots ,a_{m-1},b_1,b_2\}\cup U)$-collection in $G$ such that
$|N(X)|\le m+1$ for all $X\in\mathcal{X}_2'$. Note that
$V(G/\mathcal{X}_2')=V(G/\mathcal{X}_2)$ and
$E(\mathcal{G}/\mathcal{X}_2')\subseteq E(\mathcal{G}/\mathcal{X}_2)$. If
$E(\mathcal{G}/\mathcal{X}_2')\subsetneq E(\mathcal{G}/\mathcal{X}_2)$, then \begin{align*}
    e(\mathcal{G}/\mathcal{X}_2')\le e(\mathcal{G}/\mathcal{X}_2)-1
&\le (m+1)v(G/\mathcal{X}_2)-m^2/2-3m/2-|U|-1\\
&=(m+1)v(G/\mathcal{X}_2')-m^2/2-3m/2-|U|-1,
\end{align*}
so $\mathcal{X}_2'$ shows that Theorem
\ref{criticallylinked} holds for ${\cal G}$ and $U$, a contradiction. Therefore,
$\mathcal{G}/\mathcal{X}_2'=\mathcal{G}/\mathcal{X}_2$; so  $\mathcal{X}_2'$ satisfies (6) as well. Since $|\bigcup_{X\in\mathcal{X}_2}X|=|\bigcup_{X\in\mathcal{X}_2'}X|$ and $|\mathcal{X}_2'|>|\mathcal{X}_2|$,  $\mathcal{X}_2'$ contradicts our choice of $\mathcal{X}_2$.
\qed 

\medskip 

We claim that \begin{itemize}
    \item[(8)] $\mathcal{G}/\mathcal{X}_2-\{u_1,a_1,\ldots,a_{m-1}\}$ does not contain a $b_1$-$b_2$ path.
\end{itemize}
Suppose for a contradiction that $B^*$ is a $b_1$-$b_2$ path in
$\mathcal{G}/\mathcal{X}_2-\{u_1,a_1,\ldots,a_{m-1}\}$. We may choose
$B^*$ to be induced. Let $F$ be the set of edges $e\in E(B^*)$ such that both incident
 vertices of $e$, denoted  $l_e$ and $r_e$,  are contained in $N(X_e)$
 for some $X_e\in\mathcal{X}_2$.
 Since $B^*$ is induced and $N(X)$ induces a clique in $\mathcal{G}/\mathcal{X}_2$ for every $X\in\mathcal{X}_2$, we have $X_e\ne X_{e'}$ for any distinct $e,e'\in F$.

 For each $e\in F$ with ${N(X_e)\choose 2}\cap\{\{u_1,a_i\}:i \in [m-1]\}=
 \emptyset$, let $B_e$ be an $l_e$-$r_e$ path in $G[X_e\cup
 \{l_e,r_e\}]$ (by (7)). 
For each $e\in F$ with ${N(X_e)\choose 2}\cap\{\{u_1,a_i\}:i \in
[m-1]\}\ne \emptyset$,
let $H_e:=G[N[X_e]]$ and $t_e=|N(X_e)|\le m+1$. Note that $\mathcal{H}_e:=(H_e,
N(X_e)\setminus\{l_e,r_e\},l_e,r_e)$ is a $(t_e-2)$-rooted graph.  
If  $\mathcal{H}_e$ is feasible, let $(A_e,B_e)$ be a linkage pair in
$\mathcal{H}_e$. 

Suppose $\mathcal{H}_e$ is feasible for all  $e\in F$ with ${N(X_e)\choose 2}\cap\{\{u_1,a_i\}:i \in
[m-1]\}\ne \emptyset$. Then let  $B'$ be the path obtained from $B^*$
by, for each $e\in F$, replacing $e$ by the path $B_e$. Since
$u_1\not\in B'$, $G-B'$ contains no connected subgraph containing $a_1,\ldots,a_{m-1}$; so
there
exists some $i\in [m-1]$ such that $G-B'$ contains no $u_1$-$a_i$
path. In particular, $ua_i\notin E(G)$. Moreover,
$\{u_1,a_i\}\not\subseteq N(X)$ for any $X\in {\cal X}_2\setminus \{X_e:e\in F\}$;
for otherwise, since $G[X]$ is connected,
$G[X\cup\{u_1,a_i\}]$ contains a $u_1$-$a_i$ path. Thus,
$\{u_1,a_i\}\subseteq N(X_e)$ for some $e\in F$ (as $a_i\in N_{{\cal
    G}/{\cal X}_2}(u_1)$). But then, since $\{u_1, a_1, \ldots,
a_{m-1}\}\cap V(B^*)=\emptyset$, 
$\{u_1,a_i\}\subseteq N(X_e)\setminus\{l_e,r_e\}\subseteq V(A_e)$, and
$A_e$ contains a $u_1$-$a_i$ path. This is a contradiction as
$A_e\subseteq G-B'$.

Hence, there exists $e\in F$ such that $\mathcal{H}_e$ is not
feasible.  Since Theorem \ref{2m-link} holds for
$\mathcal{H}_e$ (by assumption), there exists some $N(X_e)$-collection $\mathcal{Y}$ in $H_e$ such that $|N(Y)|=N_{H_e}(Y)\le t_e-1<m+1$ for all $Y\in\mathcal{Y}$ and \[
e(\mathcal{H}_e/\mathcal{Y})\le (t_e-1)v(H_e/\mathcal{Y})-t_e^2/2+t_e/2.
\]
Let $\mathcal{X}_2':=(\mathcal{X}_2\setminus\{X_e\})\cup\mathcal{Y}$; then $\mathcal{X}_2'$ is a $(\{a_1, \ldots ,a_{m-1},b_1,b_2\}\cup U)$-collection in $G$ with  $|N(X)|\le m+1$ for all $X\in\mathcal{X}_2'$. 
Since $|N(Y)|\le t_e-1<|N(X_e)|$ for all $Y\in\mathcal{Y}$, we have $\mathcal{Y}\ne\{X_e\}$. Then because $G[X_e]$ is connected, $v(H_e/\mathcal{Y})\ge t_e+1$.
Note that $\mathcal{G}/\mathcal{X}_2'\subseteq \mathcal{G}/\mathcal{X}_2\cup \mathcal{H}_e/\mathcal{Y}$ and 
$\binom{N(X_e)}{2}\setminus\{l_er_e\}\subseteq E(\mathcal{G}/\mathcal{X}_2)\cap E(\mathcal{H}_e/\mathcal{Y})$.
Then
\begin{align*}
  & \quad e(\mathcal{G}/\mathcal{X}_2')\\
  &\le  e(\mathcal{G}/\mathcal{X}_2)+e(
    \mathcal{H}_e/\mathcal{Y})-\left( \binom{t_e}{2}-1\right)\\
    &\le (m+1)v(G/\mathcal{X}_2)-m^2/2-3m/2-|U|+(t_e-1)v(H_e/\mathcal{Y})-t_e^2/2+t_e/2 -\left( \binom{t_e}{2}-1\right)\\
    &= (m+1)(v(G/\mathcal{X}_2)+ v(H_e/\mathcal{Y})-t_e )-(m+2-t_e)v(H_e/\mathcal{Y})+mt_e+2t_e -t_e^2+1 \\
    &\quad -m^2/2-3m/2-|U| \\
    &\le (m+1)v(G/\mathcal{X}_2')-m^2/2-3m/2-|U|-(m+1-t_e)(v(H_e/\mathcal{Y})-t_e)-(v(H_e/\mathcal{Y})-t_e-1 )\\
    &\le (m+1)v(G/\mathcal{X}_2')-m^2/2-3m/2-|U|.
\end{align*}
If $e(\mathcal{G}/\mathcal{X}_2')<
(m+1)v(G/\mathcal{X}_2')-m^2/2-3m/2-|U|$, then $\mathcal{X}_2'$ shows
that Theorem \ref{criticallylinked} holds for ${\cal G}$ and $U$, a
contradiction. So the above inequalities are in fact equalities. In particular,\[
e(\mathcal{G}/\mathcal{X}_2')=  e(\mathcal{G}/\mathcal{X}_2)+e( \mathcal{H}_e/\mathcal{Y})-\binom{t_e}{2}+1,
\]
which implies $\binom{N(X_e)}{2}\setminus\{l_er_e\}\subseteq E(\mathcal{G}/\mathcal{X}_2')$. 
Because we also have $E(\mathcal{G}/\mathcal{X}_2)\setminus\binom{N(X_e)}{2}\subseteq E(\mathcal{G}/\mathcal{X}_2')$ and $\{u_1a_1,\ldots,u_1a_{m-1}\}\subseteq E(\mathcal{G}/\mathcal{X}_2)$,
$E(\mathcal{G}/\mathcal{X}_2')$ contains
$\{u_1a_1,\ldots,u_1a_{m-1}\}$ as well, and hence ${\cal X}_2'$ satisfies (6). However, $|\bigcup_{X\in\mathcal{X}_2'}X|<|\bigcup_{X\in\mathcal{X}_2}X|$, so $\mathcal{X}_2'$ contradicts our choice of $\mathcal{X}_2$.
\qed

\begin{itemize}
    \item[(9)]$G-\{u_1,a_1,\ldots,a_{m-1}\}$ does not contain a $b_1$-$b_2$ path.
\end{itemize}
For, if $B'$ is a $b_1$-$b_2$ path in $G-\{u_1,a_1,\ldots,a_{m-1}\}$,
then it is easy to see that in $G/\mathcal{X}_2$, $B'$ gives rise to
a $b_1$-$b_2$ path avoiding $\{u_1,a_1,\ldots,a_{m-1}\}$, which contradicts (8).\qed 

\medskip

Let $D$ be the component of $G-\{u_1,a_1,\ldots,a_{m-1}\}$ containing $b_1$. Note that $B[b_1,u_1) \subseteq D$ and $B[u_1,b_2]\subseteq G-D$.
Let $C_1,\ldots,C_t$ be the components of $D-B[b_1,u_1)$. Then ${\cal
  C}=\{V(C_i): i\in [t]\}$ is a $(\{a_1, \ldots, a_{m-1},
b_2\}\cup U)$-collection in $G-B[b_1,u_1)$ and, in fact,
$N_{G-B[b_1,u_1)}(C_i)\subseteq \{a_1, \ldots, a_{m-1},u_1\}$. Write
$G'=(G-B[b_1,u_1))/{\cal C}$, and note that $V(G')=V(G)\setminus V(D)$
and \[
E(G')\setminus E(G)\subseteq \bigcup_{j=1}^t\binom{N_{G-B[b_1,u_1)}(C_j)}{2}
\subseteq \binom{\{u_1,a_1,\ldots,a_{m-1}\}}{2}.
\]
We claim that

\begin{itemize}
    \item[(10)] the $(m-1)$ rooted graph $\mathcal{G}':=(G',\{a_1,\ldots,a_{m-1}\},u_1,b_2)$ is critically feasible with respect to $U\setminus\{u_1\}$.
\end{itemize}
It is easy to see that  $(G'[V(A)\setminus V(D)], B[u_1,b_2])$ is a linkage pair in $\mathcal{G}'$; so $\mathcal{G}'$ is feasible.
Let $(A_1,B_1)$ be an arbitrary linkage pair in ${\cal G}'$. Since $B_1$ does not contain any edge in $ \binom{\{u_1,a_1,\ldots,a_{m-1}\}}{2}$, we have $B_1\subseteq G$. Since $u_1\not\in A_1$, we have \[
E(A_1)\setminus E(G)\subseteq \{a_{i_1}a_{i_2}: \{i_1, i_2\}\in
{[m-1]\choose 2} \mbox{ and } \{a_{i_1},a_{i_2}\}\subseteq
N_{G-B[b_1,u_1)}(C_j) \mbox{ for some $j\in [t]$} \}.
\]
Let $A'$ be the graph obtained from $A_1$ by, for each $a_{i_1}a_{i_2}\in E(A_1)\setminus E(G)$, replacing the edge $a_{i_1}a_{i_2}$ in $A_1$ by an $a_{i_1}$-$a_{i_2}$ path in $G[V(C_j)\cup\{a_{i_1},a_{i_2}\}]$, where $\{a_{i_1},a_{i_2}\}\subseteq N_{G-B[b_1,u_1)}(C_j)$. Then $A'$ is a connected subgraph of $G$, and $(A',B[b_1,u_1]\cup B_1)$ form a linkage pair in $\mathcal{G}$. 
Because $\mathcal{G}$ is critically feasible with respect to $U$, $U\subseteq V(B[b_1,u_1]\cup B_1)$, which implies $U\setminus\{u_1\}\subseteq V(B_1)$. Therefore, $\mathcal{G}'$ is critically feasible with respect to $U\setminus\{u_1\}$.
\qed 

\medskip

Now by (10), we may apply Theorem \ref{criticallylinked} inductively
to $\mathcal{G}'$ and $U\setminus\{u_1\}$. Then there exists a $(\{a_1, \ldots ,a_{m-1},b_2\}\cup U)$-collection $\mathcal{X}'$ in $G'$ such that
$|N_{G'}(X)|\le m+1$ for all $X\in\mathcal{X}'$, and
\begin{align*}
   e(\mathcal{G}'/\mathcal{X}')&\le (m+1)v(G'/\mathcal{X}')-m^2/2-3m/2-1-|U\setminus\{u_1\}| \\
   &=(m+1)v(G'/\mathcal{X}')-m^2/2-3m/2-|U|.
\end{align*}
Let $\mathcal{X}=\mathcal{X}'\cup\{V(D)\setminus\{b_1\}\}$.
Then  $\mathcal{X}$
is a $(\{a_1, \ldots ,a_{m-1},b_1,b_2\}\cup U)$-collection in $G$, 
 $|N(X)|\le m+1$ for all $X\in\mathcal{X}$,  
$$
V(G/\mathcal{X})= V(G'/\mathcal{X}')\cup\{b_1\},\mbox{ and } 
E(\mathcal{G}/\mathcal{X})\subseteq   E(\mathcal{G}'/\mathcal{X}')\cup\{b_1u_1,b_1a_1,\ldots,b_1a_{m-1}\}.$$
Hence,  \begin{align*}
     e(\mathcal{G}/\mathcal{X})&\le e(\mathcal{G}'/\mathcal{X}')+m\\ 
     &\le (m+1)v(G'/\mathcal{X}')-m^2/2-3m/2-|U|+m\\
      &= (m+1)(v(G/\mathcal{X})-1)-m^2/2-3m/2-|U|+m\\
     &=(m+1)v(G/\mathcal{X})-m^2/2-3m/2-1-|U|.
 \end{align*}
 So $\mathcal{X}$ shows that Theorem
 \ref{criticallylinked} holds for ${\cal G}$ and $U$, which
 contradicts our assumption.

 \medskip

 To summarize the discussion of this section, we have  the following
 lemma. 

\begin{lemma}\label{critical}
Let $m\ge 2$ be an integer. Suppose Theorem \ref{2m-link} holds for all $m'$-rooted
graphs with $0\le m'\le m-1$, and suppose Theorem
\ref{criticallylinked} holds for all critically feasible
$(m-2)$-rooted graphs.
Then Theorem \ref{criticallylinked} holds for all critically feasible $(m-1)$-rooted graphs.
\end{lemma}

\section{Conclusion}

First, we complete the proofs of Theorem \ref{2m-link} and Theorem
\ref{criticallylinked}. We apply induction on $m\ge 0$. By Lemma
\ref{base1}, Theorem \ref{2m-link} holds for $m=0,1$ or $2$. By
Lemma~\ref{base2},  Theorem \ref{criticallylinked} holds for $m=0$;
and by Lemma \ref{critical}, Theorem \ref{criticallylinked} holds for $m=1$ and $m=2.$

Now assume $m\ge 3$ and suppose Theorem \ref{2m-link} 
and Theorem \ref{criticallylinked} both
hold for  $m'$-rooted graphs with $0\le m'\le m-1$. 
By Lemma \ref{2m}, Theorem \ref{2m-link} holds for all $m$-rooted graphs. Hence, by Lemma \ref{critical}, Theorem \ref{criticallylinked} holds for all critically feasible $m$-rooted graphs. \qed 

\medskip

Next, we use Theorem~\ref{2m-link} to prove Theorem~\ref{main}. Let $G$ be a $(2m+2)$-connected graph, and let $a_1,\ldots, a_m,b_1,b_2$ be distinct vertices in $G$. By Theorem~\ref{2m-link}, $G$ contains a $b_1$-$b_2$ path $B$ such that $\{a_1, \ldots, a_m\}$ is contained in a component of $G-B$, which we denote by $A(B)$. Let $C_1,\ldots, C_t$ denote the components of $G-B$ such that $C_1=A(B)$ and $|V(C_2)|\ge |V(C_3)|\ge \ldots \ge |V(C_t)|$. We choose $B$ such that $(|V(C_1)|, |V(C_2)|, \ldots, |V(C_t)|)$ is maximal with respect to the lexicographic ordering.

We claim  $t=1$ and, hence, $G-B$ is connected. For, suppose $t\ge
2$. Let $u_1,u_2\in N(C_t)\cap V(B)$ with $B[u_1,u_2]$ maximal. Since
$G$ is $(2m+2)$-connected, it is at least 3-connected; so there exist
$u\in V(B(u_1,u_2))$ and $s\in [t-1]$ such that $N(u)\cap C_s\ne
\emptyset$. We choose $u,s$ so that $s$ is minimum. Let $B'$ be
obtained from $B$ by deleting $B(u_1,u_2)$ and adding an induced
$u_1$-$u_2$ path in $G[V(C_t)\cup \{u_1,u_2\}]$. Now $B'$ is a
$b_1$-$b_2$ path in $G$ and $C_1$ is contained in a component of
$G-B'$. If we let $C_1',\ldots, C_q'$ denote the components of $G-B'$
such that $C_1\subseteq C_1'$ and $|V(C_2')|\ge |V(C_3')|\ge \ldots
\ge |V(C_{q})'|$, then 
$|V(C_i')|\ge |V(C_i)|$ for $i\in [s-1]$ and $|V(C_s')|>|V(C_s)|$. This contradicts the choice of $B$ that $(|V(C_1)|, |V(C_2)|, \ldots, |V(C_t)|)$ is maximal, completing the proof of Theorem~\ref{2m-link}. \qed

\medskip

Question~\ref{removable} is related to the following conjecture in \cite{HKY15}.

\begin{conjecture}[Hong, Kang, Yu] \label{hky}
  For any positive integers $k$ and $l$, there exists a smallest positive integer $g(k, l)$
such that for any $g(k, l)$-connected graph $G$, any edge $e \in E(G)$, and any $k$-element set $X \subseteq V(G)$ that is not incident with $e$, there exists an induced cycle $C$ in $G -X$ such that $e\in E(C)$ and $G-C$ is
$l$-connected.
\end{conjecture}

Conjecture~\ref{hky} for $k=0$ is an old conjecture made by Lov\'{a}sz
\cite{Lo75} in 1975. The existence of $g(0,l)$ for $l\ge 3$  is still
open, while $g(0,1)=3$ follows from a result of Tutte \cite{Tu63} and $g(0,2)=5$ is proved independently in \cite{CGY03} and \cite{Kr01}.

It is shown in \cite{HKY15} that $g(k,1) \le 10k+1$. Theorem~\ref{main} implies that $g(k,1)\le 2k+3$. It would be interesting to also improve the bound $g(k,2)\le 10k+11$ in \cite{HKY15}.

\end{document}